# A Lagrangian Decomposition Algorithm for Robust Green Transportation Location Problem

Ali Rouhani, Mahdi Bashiri*, Rashed Sahraeian

*bashiri@shahed.ac.ir

Department of Industrial Engineering, Shahed University, Tehran, Iran

**Abstract:**

In this paper, a green transportation location problem is considered with uncertain demand parameter. Increasing robustness influences the number of trucks for sending goods and products, and consequently, makes the air pollution enhance. In this paper, two green approaches are introduced which demand is the main uncertain parameter in both. These approaches are addressed to provide a trade-off between using available trucks and buying new hybrid trucks for evaluating total costs beside air pollution. Due to growing complexity, a Lagrangian decomposition algorithm is applied to find a tight lower bound for each approach. In this propounded algorithm, the main model is decomposed into master and subproblems to speed up convergence with a tight gap. Finally, the suggested algorithm is compared with commercial solver regarding total cost and computational time. Due to computational results for the proposed approach, the Lagrangian decomposition algorithm is provided a close lower bound in less time against commercial solver.

## 1. INTRODUCTION

Increasing the number of required products and developing transportation systems are the main result of population growth. As a result, it makes the air polluted, and mechanisms for controlling pollution become important. When companies produce a particular product in the first time, have not any accessible data about the products' demand. Therefore, they should estimate their volume. There are several ways to dealing with this uncertainty. In recent years, one of the ways that have progressed remarkably is robust optimization. Robust optimization generally divides into two types of interval-based and scenario-based models. In this paper, interval-based robust optimization is considered. In this term, there are some pioneers such as Soyster [1], Ben-Tal [2], and Bertsimas [3] approaches. Bertsimas's approach used regarding its flexibility on considering uncertain parameters related to other approaches [4, 5, 6]. Budget parameter effects on the price of robustness in this approach. The number of trucks and consequently pollution emission is enhanced by increasing the amount of the budget parameter. Some companies eliminate these problems by choosing the costly solution and buy new hybrid trucks. However, controlling pollution created by available trucks is an economical solution against the first approach. These two approaches are examined in this paper. For considering the second approach, suppose that published pollution of each truck is followed from a distribution function with known mean and variance, and can be controlled with a threshold. Under this assumption, a chance constraint is necessary to be used. The chance constraint is one of the hard and probabilistic restrictions, which can be added to the main problem.

In large-scale MIP problems, commercial solvers' efficiency is reduced. Therefore, decomposition algorithms may be used to solve these problems. Decomposition-based solution methods are employed to find exact solutions for MIP problems. In the contrast of other decomposition-based algorithms, Lagrangian decomposition algorithm is considered to find a tight lower bound for large-scale problems. Lagrangian decomposition is a kind of Lagrangian relaxation algorithm which decomposes the problem into some subproblems after relaxing hard constraints. For using Lagrangian decomposition, some methods have been introduced previously like the subgradient method, cutting plane and so on which, in this paper, the second one method is applied. In this conception, the master problem is considered for reducing iteration of solving problems after decomposing the model into two subproblems. In the master problem, the main decision variables are fixed, and Lagrange multipliers are found as decision variables.

The rest of this paper is organized as follows: In the next section, a literature review is presented for green transportation location problem with uncertain demands. In section 3, the proposed mathematical model is presented. In the fourth section, the

Lagrangian decomposition and steps of this algorithm are discussed. Sensitivity analysis and computational experiments are examined in section 5. Finally, the paper concludes in the last section.

## 2. LITERATURE REVIEW

In this section, the published paper are reviewed and contribution of paper is presented according to its evaluation. A two stage robust was mentioned by [7]. A transportation location problem has modeled with two stage stochastic programming concept, which distribution channels are a priori decisions to optimized network flow. Second stage variables are flow and origins decisions [8]. A novel mathematical model for transportation location problem has presented in disaster application, which location and origin-destination allocation decisions are priori known [9]. A two stage stochastic programming is used to deal with parameter uncertainty. First stage variables are flow of priori allocation and new allocation decisions. Flow of new distribution channels and shortage or leftovers of distribution channels are second stage variables [10]. A bi objective mixed integer location/routing model have presented that aims to minimize transportation cost and risks for large-scale hazardous waste management systems (HWMSs), whereas all parameters are known [11]. A summary of reviewed papers are compared in table 1. Main contributions of this study can be summarized as follows;

- Controlling amount of pollution in the network with a chance constraint concept.
- Using robust optimization for dealing with uncertainty in the green transportation location problem.
- Considering a Lagrangian decomposition algorithm for the robust green transportation-location problem.

| TABLE 1. A summary of related literature review ||||||||
|---|---|---|---|---|---|---|---|
| No. | Reference | Year | App. | D. or U. | Unc. Form. | No. Obj | Solution algorithm |
| 1 | [9] | 2017 | Network Flow | D | - | 1 | SA |
| 2 | [10] | 2015 | Disaster | D | - | 1 | CS |
| 3 | [11] | 2016 | Hazardous waste | D | - | 2 | CS |
| 4 | [12] | 2014 | Disaster | D | - | 1 | CS |
| 5 | [8] | 2014 | Network Flow | U | TSSP | 1 | CS |
| 6 | [7] | 2014 | Network Flow | U | TSRO | 1 | Kelley's Algorithm |
| 7 | [13] | 2004 | Disaster | U | TSSP | 1 | CS |
| 8 | [14] | 2009 | Network Flow | U | TSSP | 1 | Monte Carlo |
| 9 | [15] | 2015 | Disaster | U | TSSP | 1 | GA |
| 10 | [16] | 2014 | Network Flow | D | - | 2 | NSGA II |
| 11 | [17] | 2014 | Disaster | D | - | 3 | ECA |
| 12 | [18] | 2011 | Disaster | D | - | 2 | NSGA II |
| App: Application, D or U: Deterministic or Uncertain, Unc Form: types of uncertainty formulation, No Obj: Number of Objectives, CS: Commercial Solver, TSSP: Two stage stochastic programming, TSRO: Two stage Robust Optimization ||||||||

Also, Lagrangean decomposition has been used in various problems such as Quadratic binary Program [20], location-allocation problem, offshore oilfield development planning [21] as a solution algorithm. Due to mentioned papers, we apply Lagrangean decomposition for this problem which not used until now.

## 3. PROBLEM DEFINITION

This section is divided into two parts, in the first part, a transportation location problem is defined with demand uncertainty. In the second one, green approaches are mentioned.

A transportation-location problem is composed of transportation and location-allocation problems, and its aim is transporting each product due to the amount of demand in each destination with the minimum total cost. The capacity of origins and trucks restrict sending products, and it is assumed that the required

vehicles already exist in the shipping company. Considered problem costs are included:

A. Shipping costs from the origins to destinations

B. The cost of linking between origins and destinations

C. The cost of established origins

D. The cost of the shortage of products at destinations

For example, a company is planned to produce various products and deliver to customers with regard to the total cost. The company should use different types of trucks for satisfying customers' demands. The capacity of the origins and trucks are playing an important role in the number of delivered products.

Due to the mentioned example, suppose that this company produces new commodities; while, market demands are unknown, and company revenue is increased when all market demands are satisfied. If the company wants to satisfy all customer demands, the number of trucks and consequently, pollution emissions are increased. Two approaches are suggested for addressing pollution emissions: Firstly, due to required trucks, the company is decided to purchase hybrid trucks for sending products to destinations. In the second approach, the pollution emission is controlled by adding some limitations regarding as the age of trucks. Assumed that it follows a normal distribution, so related constraints are added as the chance constraint.

## 4. MATHEMATICAL MODEL

**Sets**

$I$: The set of Origins
$J$: The set of destinations
$P$: The set of trucks
$L$: The set of products

**Parameters**

$c_{ij}^p$: Set up cost for link between $i$ th origin to $j$ th destination for truck $p$

$q_{ij}^{lp}$: Transportation cost for $l$ th product with $p$ th truck from $i$ th origin to $j$ th destination

$h_i^l$: $i$ th origin opening cost for $l$ th product

$w_j^l$: Penalty cost for unmet demands for $l$ th product in $j$ th destination

$b^{lp}$: Capacity of $p$ th truck for $l$ th product

$k_i^l$: Capacity of $i$ th origin for $l$ th product
$em^p$: Pollution emission from $p$ th truck
$cbc^p$: Purchasing cost of $p$ th hybrid truck
$D_j^l$: Demands of $l$ th product in $j$ th destination
$Td_{ij}$: Maximum allowable pollution emission in link between $i$ and $j$
$1-\alpha$: Confidence interval

**Variables**

$y_{ij}^p$: 1 If a link between $i$ and $j$ is constructed for truck $p$ and 0, otherwise

$x_{ij}^{lp}$: Flow between $i$ and $j$ by the truck $p$ for $l$ th product

$z_i^l$: 1 if origin $i$ is used for shipping commodity $l$, 0 otherwise

$u_j^l$: Amount of unsatisfied demand in $j$ th destination for $l$ th product

$nc^p$: Number of $p$ th needed hybrid trucks
$ofv1$: Objective function 1 (Total Cost)
$ofv2$: Objective function 2 (Total Cost)

### 4.1. Mathematical problem for the first approach:

$$Min\ ofv1 = \sum_{i=1}^{I}\sum_{j=1}^{J}\sum_{p=1}^{P} c_{ij}^p y_{ij}^p +$$

$$\sum_{i=1}^{I}\sum_{l=1}^{L} h_i^l z_i^l + \sum_{i=1}^{I}\sum_{j=1}^{J}\sum_{p=1}^{P}\sum_{l=1}^{L} q_{ij}^{lp} x_{ij}^{lp} + \quad (1)$$

$$\sum_{j=1}^{J}\sum_{l=1}^{L} w_j^l u_j^l + \sum_{p=1}^{P} cbc^p nc^p$$

$$s.t.\ \sum_{i=1}^{I}\sum_{p=1}^{P} y_{ij}^p \geq 1,\quad \forall j \quad (2)$$

$$\sum_{l=1}^{L} x_{ij}^{lp} \leq \sum_{l=1}^{L} b^{lp} y_{ij}^p,\quad \forall i,j,p \quad (3)$$

$$\sum_{j=1}^{J}\sum_{p=1}^{P} x_{ij}^{lp} \leq k_i^l z_i^l,\quad \forall i,l \quad (4)$$

$$u_j^l \geq D_j^l - \sum_{i=1}^{I}\sum_{p=1}^{P} x_{ij}^{lp},\quad \forall j,l \quad (5)$$

$$nc^p = \sum_{i=1}^{I}\sum_{j=1}^{J} y_{ij}^p,\quad \forall p \quad (6)$$

$$y_{ij}, z_i^l \in \{0,1\}$$
$$x_{ij}^{lp}, u_j^l, nc^p \geq 0 \quad (7)$$

## 4. 2. Mathematical problem for the second approach:

$$Min\ ofv2 = \sum_{i=1}^{I}\sum_{j=1}^{J}\sum_{p=1}^{P} c_{ij}^{p} y_{ij}^{p} +$$

$$\sum_{i=1}^{I}\sum_{l=1}^{L} h_{i}^{l} z_{i}^{l} + \sum_{i=1}^{I}\sum_{j=1}^{J}\sum_{p=1}^{P}\sum_{l=1}^{L} q_{ij}^{lp} x_{ij}^{lp} \quad (8)$$

$$+ \sum_{j=1}^{J}\sum_{l=1}^{L} w_{j}^{l} u_{j}^{l}$$

$$s.t.\ \sum_{i=1}^{I}\sum_{p=1}^{P} y_{ij}^{p} \geq 1, \quad \forall j \quad (9)$$

$$\sum_{l=1}^{L} x_{ij}^{lp} \leq \sum_{l=1}^{L} b^{lp} y_{ij}^{p}, \quad \forall i,j,p \quad (10)$$

$$\sum_{j=1}^{J}\sum_{p=1}^{P} x_{ij}^{lp} \leq k_{i}^{l} z_{i}^{l}, \quad \forall i,l \quad (11)$$

$$u_{j}^{l} \geq D_{j}^{l} - \sum_{i=1}^{I}\sum_{p=1}^{P} x_{ij}^{lp}, \quad \forall j,l \quad (12)$$

$$\Pr(\sum_{p=1}^{P} em^{p} y_{ij}^{p} \leq Td_{ij}) \geq 1-\alpha, \quad \forall i,j \quad (13)$$

$$y_{ij}, z_{i}^{l} \in \{0,1\}$$
$$x_{ij}^{lp}, u_{j}^{l} \geq 0 \quad (14)$$

For clarifying the mathematical problem, note that some equations are duplicated in each formulation. These equations are mentioned in the pair number which the first one is about the first approach formulation and the second one is about the second approach formulation. In the revised version, the definition of indices, parameters, variables, and equations are clarified. For example, the definition of equations was modified as following:

Equations (1) and (8) calculate the total cost of the transportation system that the two first parts are about the cost of establishing origins and destinations. The third part calculates the transportation cost between origin and destination, and the last part calculates penalty cost of unmet demands in destinations. In the first approach mathematical formulation, the capacity limitation for trucks and origins are mentioned in equations (3) and (4); however, equations (10) and (11) consider the second approach capacity restrictions for origin and trucks. Amount of unsatisfied demands are determined in equations (5) and (12). Equation (6) calculates the number of each truck which was used for the first conception. In the conception, pollution emission threshold is considered in the equation (13).

## 4. 3. Robust optimization in green transportation location problem:

Bertsimas and Sim approach [3] can control effects of demand uncertainty on the network design. This approach was presented by Bertsimas et al [3] and improved [19]. In the improved paper, a novel approach is used to model the robust green closed loop supply chain problem [19].

Assume that demand interval is $[D_{j}^{l} + \Delta D_{j}^{l-}, D_{j}^{l} + \Delta D_{j}^{l+}]$ that $D_{j}^{l}$ it is the nominal value of demand and $\Delta D_{j}^{l+}, \Delta D_{j}^{l-}$ are positive and negative deviation from the nominal value, respectively. Assuming that computing unmet demand penalty cost is computed in the following constraint:

$$\sum_{j=1}^{J}\sum_{l=1}^{L} w_{j}^{l}(D_{j}^{l} - \sum_{i=1}^{I}\sum_{p=1}^{P} x_{ij}^{lp}) \leq v \quad (15) \quad (16)$$

$v$ is a free variable which used in the objective function instead of $\sum_{j=1}^{J}\sum_{l=1}^{L} w_{j}^{l} u_{j}^{l}$.

Robust counterpart constraints are presented as bellow:

$$\sum_{j=1}^{J}\sum_{l=1}^{L} w_{j}^{l}(D_{j}^{l} - \sum_{i=1}^{I}\sum_{p=1}^{P} x_{ij}^{lp}) + \underset{\beta_{j}^{l+},\beta_{j}^{l-} \in k}{Max}[$$
$$\sum_{j=1}^{J}\sum_{l=1}^{L} (\Delta D_{j}^{l+}\beta_{j}^{l+} - \Delta D_{j}^{l-}\beta_{j}^{l-}) w_{j}^{l}] \leq v \quad (16)$$

According to duality theorem, nonlinear constraints can be rewritten regarding to:

1. $w_{j}^{l}$ Constant is not necessary to use in dual transform.

2. $\beta_{j}^{l+}$ And $\beta_{j}^{l-}$ are coefficients of positive and negative deviation that in the worst case may be 1.

Linearized counterparts are as followed;

$$- \underset{\beta_j^{l+}, \beta_j^{l-} \geq 0 \, \& \in k}{Min} [\sum_{j=1}^{J}\sum_{l=1}^{L}(-\Delta D_j^{l+}\beta_j^{l+} + \Delta D_j^{l-}\beta_j^{l-})] \quad (17)$$

$$-\beta_j^{l+} \geq -1, \quad \forall J, L \quad (\alpha 1_j^l) \quad (18)$$

$$-\beta_j^{l-} \geq -1, \quad \forall J, L \quad (\alpha 2_j^l) \quad (19)$$

$$-\beta_j^{l+} - \beta_j^{l-} \geq -\Gamma_j^l, \quad \forall J, L \quad (\mu_j^l) \quad (20)$$

Dual variables are shown in the parentheses. Finally dual model are mentioned bellow:

$$\underset{\alpha 1_j^l, \alpha 2_j^l, \mu_j^l \geq 0}{Max} [\sum_{j=1}^{J}\sum_{l=1}^{L}(-\alpha 1_j^l - \alpha 2_j^l - \Gamma_j^l \mu_j^l)] \quad (21)$$

$$-\alpha 1_j^l - \mu_j^l \leq -\Delta D_j^{l+}, \quad \forall J, L \quad (22)$$

$$-\alpha 2_j^l - \mu_j^l \leq \Delta D_j^{l-}, \quad \forall J, L \quad (23)$$

Final robust counterparts are written as follows:

$$\sum_{j=1}^{J}\sum_{l=1}^{L}((D_j^l + \alpha 1_j^l + \Gamma_j^l \mu_j^l - \sum_{i=1}^{I}\sum_{p=1}^{P}x_{ij}^{lp})w_j^l) \leq v \quad (24)$$

$$\alpha 1_j^l + \mu_j^l \geq \Delta D_j^{l+}, \quad \forall J, L \quad (25)$$

$$\alpha 2_j^l + \mu_j^l \geq \Delta D_j^{l-}, \quad \forall J, L \quad (26)$$

$$\alpha 1_j^l, \alpha 2_j^l, \mu_j^l, v \geq 0 \quad (27)$$

### 4. 4. Chance Constraint Programming in green transportation location problem:

Pollution emission constraint is propounded by a chance constraint format, the linearization of these constraints is considered as below:

Suppose that pollution emission of each truck ($em^p$) has a normal distribution with ($E(em^p)$, $VAR(em^p)$). The linearized constraints are given in equation (29):

$$\sum_{p=1}^{P}E(em^p)y_{ij}^p + Z_{1-\alpha}\sqrt{\sum_{p=1}^{P}VAR(em^p)(y_{ij}^p)^2} \leq Td_{ij}, \quad \forall i, j \quad (28)$$

### 5. LAGRANGIAN DECOMPOSITION

Used Lagrangian decomposition is based on [20], which considered cutting planes to provide a tight lower bound. In this section, both component of Lagrangian decomposition and pseudo-code of the algorithm is presented, respectively. For starting this algorithm, relaxed constraint must be determined. Relaxed constraint is:

$$\sum_{l=1}^{L}x_{ij}^{lp} \leq \sum_{l=1}^{L}b^{lp}y_{ij}^p, \quad \forall i, j, p \quad (29)$$

After relaxing mentioned constraint, the main problem decomposed into two sub-problems that presented in the rest of the paper.

Extra parameters that used in this algorithm are mentioned below:

$\theta$ : Upper bound of first sub-problem

$\eta$ : Upper bound of second sub-problem

$\lambda_{ij}^p$ : Lagrange multipliers for relaxed constraints

### 5. 1. Lagrangean Sub Problems:

With relaxing constraint (29), Lagrangean subproblems for both approaches can be demonstrated from their main mathematical models. After relaxing mentioned constraint with Lagrange multipliers, the relaxed problem is divided into two independent problems.

### 5.2. Lagrangean Master Problems:

**First approach master problem**

$$Max \, LM = \theta + \eta \quad (30)$$

$$\theta \leq \sum_{i=1}^{I}\sum_{j=1}^{J}\sum_{p=1}^{P}(c_{ij}^p - \lambda_{ij}^p \sum_{l=1}^{L}b^{lp})y_{ij}^p + \sum_{p=1}^{P}cbc^p nc^p \quad (31)$$

$$\eta \leq \sum_{j=1}^{J}\sum_{i=1}^{I}\sum_{l=1}^{L}\sum_{p=1}^{P}(q_{ij}^{lp} + \lambda_{ij}^p)x_{ij}^{lp} + \sum_{i=1}^{I}\sum_{l=1}^{L}h_i^l z_i^l + \hat{v} \quad (32)$$

$$\lambda_{ij}^p \geq 0 \quad (33)$$

**Second approach master problem**

$$Max \, LM = \theta + \eta \quad (34)$$

$$\theta \leq \sum_{i=1}^{I}\sum_{j=1}^{J}\sum_{p=1}^{P}(c_{ij}^p - \lambda_{ij}^p \sum_{l=1}^{L} b^{lp})y_{ij}^p \qquad (35)$$

$$\eta \leq \sum_{j=1}^{J}\sum_{i=1}^{I}\sum_{l=1}^{L}\sum_{p=1}^{P}(q_{ij}^{lp} + \lambda_{ij}^p)x_{ij}^{lp} + \sum_{i=1}^{I}\sum_{l=1}^{L} h_i^l z_i^l + \hat{v} \qquad (36)$$

$$\lambda_{ij}^p \geq 0 \qquad (37)$$

Suggested algorithm Pseudo-code is presented as follow:

### 5.3. Lagrangian Decomposition algorithm Pseudo-code:

**1) Initialized:**
$$Z_{up} = \infty, \quad Z_{lb} = -\infty, \quad iter = 1$$

**2) Solve Lagrangian Sub problems:**
Store all variables
Store Objective functions values
If sum of the objective values are greater than $Z_{lb}$, update $Z_{lb}$

**3) Solve Lagrangian Master Problem:**
Store $\lambda_{ij}^p$
If Master problem objective function are lower than $Z_{up}$, update $Z_{up}$

**4) Convergence test:**
If $Z_{up} - Z_{lb} < \varepsilon$ stop the algorithm
Else go to step 2

## 6. NUMERICAL STUDIES

Two approaches propound for dealing with published pollution, which in the first approach, the decision maker must buy new trucks to serves the customers. But, in the second approach, the decision maker tries to design supply network somehow that the available truck published pollution is not exceeding form the particular threshold. In this section, the algorithm-based results illustrate. In the other words, the Lagrangean decomposition applies for each proper conception and computational and GAP percent of this algorithm are compared.

Algorithm-based results are illustrated in tables 2-5, which in these tables, computational time and gap of lagrangian decomposition algorithm are calculated. As is clear from the results, mentioned algorithm are provide closed lower bound for proposed model with lower running time, that obtained results are more specifically due to growing size of the problem. GAP measure can be computed as follow:

$$\%GAP = (\frac{Z_{M\,Problem} - Z_{P\,Algorithm}}{Z_{M\,Problem}})*100 \qquad (38)$$

**TABLE 2:** First size comparison computational time for first approach

|  | Main Problem | | Proposed Algorithm | | |
|---|---|---|---|---|---|
| NO. | Obj | T(s) | Obj | T(s) | % Gap |
| 1 | 7258546.335 | 0.5 | 7252045.364 | 0.173 | 0.089562985 |
| 2 | 7547085.984 | 0.592 | 7540585.014 | 0.109 | 0.086138819 |
| 3 | 7835625.634 | 0.484 | 7829124.663 | 0.156 | 0.082966837 |
| 4 | 8124165.283 | 0.608 | 8117664.313 | 0.141 | 0.080020168 |
| 5 | 8412704.933 | 0.296 | 8406203.962 | 0.14 | 0.07727563 |
| 6 | 8701244.582 | 0.281 | 8694743.612 | 0.155 | 0.074713114 |
| 7 | 8989784.232 | 0.265 | 8983283.261 | 0.14 | 0.072315092 |
| 8 | 9278323.881 | 0.281 | 9271822.911 | 0.14 | 0.070066219 |
| 9 | 9566863.531 | 0.281 | 9560362.56 | 0.171 | 0.067953 |
| 10 | 9855403.18 | 0.265 | 9848902.21 | 0.11 | 0.065963519 |
| average |  | 0.3853 |  | 0.1435 | 0.076698 |

**TABLE 3:** First size comparison computational time for second approach

| $Z_{1-\alpha} = -3$ | $Z_{1-\alpha} = 3$ |
|---|---|

| Main Problem | | Proposed Algorithm | | | Main Problem | | Proposed Algorithm | | |
|---|---|---|---|---|---|---|---|---|---|
| Obj | T(s) | Obj | T(s) | % Gap | Obj | T(s) | Obj | T(s) | % Gap |
| 7275656.23 | 25.65 | 7271088.77 | 0.79 | 0.0628 | 7275727.26 | 23.30 | 7258546.33 | 0.53 | 0.2361 |
| 7564266.91 | 23.30 | 7559634.28 | 0.62 | 0.0612 | 7564266.91 | 23.88 | 7547085.98 | 0.54 | 0.2271 |
| 7852806.56 | 23.10 | 7848173.93 | 0.62 | 0.059 | 7852806.56 | 23.12 | 7835625.63 | 0.57 | 0.2188 |
| 8141346.21 | 24.78 | 8136713.58 | 0.62 | 0.0569 | 8141346.21 | 24.24 | 8124165.28 | 0.56 | 0.211 |
| 8429885.86 | 25.42 | 8425253.23 | 0.60 | 0.055 | 8429885.86 | 24.55 | 8412704.93 | 0.54 | 0.2038 |
| 8718425.51 | 24.00 | 8713792.88 | 0.60 | 0.0531 | 8718425.51 | 23.32 | 8701244.58 | 0.54 | 0.1971 |
| 9006965.16 | 23.18 | 9002332.53 | 0.59 | 0.0514 | 9006965.16 | 23.04 | 8989784.23 | 0.56 | 0.1908 |
| 9295504.80 | 23.10 | 9290872.18 | 0.59 | 0.0498 | 9295504.80 | 23.15 | 9278323.88 | 0.56 | 0.1848 |
| 9584044.45 | 24.10 | 9579411.83 | 0.62 | 0.0483 | 9584044.45 | 25.19 | 9566863.53 | 0.56 | 0.1793 |
| 9872584.10 | 24.133 | 9867951.48 | 0.577 | 0.0469 | 9872584.10 | 23.431 | 9855403.1 | 0.57 | 0.174 |
| average | 24.076 | | 0.622 | 0.0544 | | 23.722 | | 0.553 | 0.2022 |

TABLE 4: Second size comparison computational time for first approach

| | Main Problem | | Proposed Algorithm | | |
|---|---|---|---|---|---|
| NO. | Obj | T(s) | Obj | T(s) | % Gap |
| 1 | 17108168.97 | 0.312 | 17081279.72 | 0.234 | 0.15717198 |
| 2 | 17880707.73 | 0.67 | 17853818.48 | 0.172 | 0.15038134 |
| 3 | 18653246.49 | 0.671 | 18626357.25 | 0.187 | 0.14415318 |
| 4 | 19425785.26 | 0.452 | 19398896.01 | 0.048 | 0.138420391 |
| 5 | 20198324.02 | 0.687 | 20171434.78 | 0.203 | 0.133126134 |
| 6 | 20970862.79 | 0.452 | 20943973.54 | 0.218 | 0.128221944 |
| 7 | 21743401.55 | 0.593 | 21716512.3 | 0.204 | 0.123666244 |
| 8 | 22515940.32 | 0.687 | 22489051.07 | 0.187 | 0.119423162 |
| 9 | 23288479.08 | 0.437 | 23261589.83 | 0.126 | 0.115461589 |
| 10 | 24061017.85 | 0.328 | 24034128.6 | 0.218 | 0.111754407 |
| average | | 0.529 | | 0.179 | 0.132178 |

TABLE 5: Second size comparison computational time for Second approach

| $Z_{1-\alpha} = -3$ | | | | | $Z_{1-\alpha} = 3$ | | | | |
|---|---|---|---|---|---|---|---|---|---|
| Main Problem | | Proposed Algorithm | | | Main Problem | | Proposed Algorithm | | |
| Obj | T(s) | Obj | T(s) | % Gap | Obj | T(s) | Obj | T(s) | % Gap |
| 12691067 | 135.78 | 12690733 | 0.624 | 0.002635 | 12691067 | 4.27 | 12690733 | 0.779 | 0.002635 |
| 13222631 | 138.73 | 13222296 | 0.968 | 0.002529 | 13222631 | 4.13 | 13222296 | 0.733 | 0.002529 |
| 13754195 | 103.88 | 13753860 | 0.967 | 0.002431 | 13754195 | 4.21 | 13753860 | 0.733 | 0.002431 |
| 14285758 | 51.12 | 14285424 | 0.936 | 0.002341 | 14285758 | 4.3 | 14285424 | 0.749 | 0.002341 |
| 14817322 | 201.05 | 14816988 | 0.982 | 0.002257 | 14817322 | 4.43 | 14816988 | 0.78 | 0.002257 |
| 15348886 | 43.13 | 15348552 | 1.061 | 0.002178 | 15348886 | 4.3 | 15348552 | 0.718 | 0.002178 |
| 15880450 | 71.83 | 15880116 | 0.936 | 0.002105 | 15880450 | 6.36 | 15880116 | 0.779 | 0.002105 |
| 16412014 | 44.63 | 16411679 | 0.967 | 0.002037 | 16412014 | 4.25 | 16411679 | 0.765 | 0.002037 |
| 16943578 | 723.29 | 16943243 | 0.999 | 0.001973 | 16943578 | 5.92 | 16943243 | 0.749 | 0.001973 |
| 17475141 | 1041.46 | 17474807 | 1.059 | 0.001913 | 17475141 | 4.13 | 17474807 | 0.78 | 0.001913 |
| average | 255.49 | | 0.95 | 0.00224 | | 4.63 | | 0.756 | 0.00224 |

Referred by results in above tables, when the sample size is small not only the computational time between each method is closed together, but also the GAP of the Lagrangian decomposition is less than large-scale problem. However, in the small size, the distinction between the two methods according to the time is so hard, but in the large one, the differences are seen obviously. Moreover, provided tables are demonstrated that in all instances Lagrangian decomposition's GAP is less than 1% which the suitability of this algorithm for solving such problems is shown.

These tables are clarified Lagrangean decomposition can employ in order to evaluating a lower bound for these mathematical formulations.

## 7. CONCLUSION

In this paper, a robust green transportation location problem is suggested with uncertain demands. Bertsimas methodology is used for dealing with demand uncertainty. Two approaches are mentioned for controlling pollution emissions. In the first approach, total purchasing cost of new hybrid trucks is examined. In the second one, chance constraints are added to control pollution emission by available trucks. According to numerical examples, a trade-off is performed and it is demonstrated that which one has a lower total cost. Lagrangian decomposition is presented for providing a tight lower bound in a rational time. Computational results confirm that the presented algorithm is efficient besides of low optimally gap. For future research, a conditional value-at-risk instead of robust optimization can be used due to the problem concept. In the problem, a definition can propound multi-period or dynamic system that can adjust published pollution in different time periods. Using exact algorithms such as bender's decomposition, Dantzig-Wolfe decomposition, and so on can reduce computational burden besides solving problems in the exact forms. Moreover, there are some applications such as telecommunication, electricity distribution systems, and production planning which can employ the proposed model to improve the performance of their optimization problems. Employing this formulation is suggested as a future study.